\documentclass [12pt]{amsart}
\usepackage{amsmath}
\usepackage{amssymb}

\newtheorem{thm}{Theorem}[section]
\newtheorem{cor}[thm]{Corollary}
\newtheorem{lem}[thm]{Lemma}

\theoremstyle{definition}

\theoremstyle{remark}

\numberwithin{equation}{section}

\newcommand{\beas}{\begin{eqnarray*}}
\newcommand{\eeas}{\end{eqnarray*}}
\newcommand{\bes} {\begin{equation*}}
\newcommand{\ees} {\end{equation*}}
\newcommand{\be} {\begin{equation}}
\newcommand{\ee} {\end{equation}}
\newcommand{\bea} {\begin{eqnarray}}
\newcommand{\eea} {\end{eqnarray}}
\newcommand{\ra} {\rightarrow}
\newcommand{\txt} {\textmd}
\newcommand{\ds} {\displaystyle}

\begin{document}

\title[\tiny{Segal-Bargmann transform and Paley-Wiener theorems on motion groups}]
 {Segal-Bargmann transform and Paley-Wiener theorems on motion groups}

\author{\tiny{Suparna Sen}}

\address{Department of Mathematics, Indian Institute of Science, Bangalore - 560012, India.}

\email{suparna@math.iisc.ernet.in.}

\thanks{The author was supported by Shyama Prasad Mukherjee Fellowship from Council of Scientific and Industrial Research, India.}


\begin{abstract}
We study the Segal-Bargmann transform on a motion group
$\mathbb{R}^n \ltimes K,$ where $K$ is a compact subgroup of
$SO(n).$ A characterization of the Poisson integrals associated to
the Laplacian on $\mathbb{R}^n \ltimes K$ is given. We also
establish a Paley-Wiener type theorem using the complexified
representations.

\vspace*{0.1in}

\begin{flushleft}
\vspace*{0.1in}
\end{flushleft}

\end{abstract}

\maketitle

\section{Introduction}

The Segal-Bargmann transform, also called the coherent state
transform, was developed independently in the early 1960's by
Segal in the infinite-dimensional context of scalar quantum field
theories and by Bargmann in the finite-dimensional context of
quantum mechanics on $\mathbb{R}^n.$ We consider the following
equivalent form of Bargmann's original result.

A function $f \in L^2(\mathbb{R}^n)$ admits a factorization $f(x)
= g*p_t(x)$ where $g \in L^2(\mathbb{R}^n)$ and $\ds{p_t(x) =
\frac{1}{(4\pi t)^{\frac{n}{2}}} e^{\frac{-|x|^2}{4t}}}$ (the heat
kernel on $\mathbb{R}^n$) if and only if $f$ extends as an entire
function to $\mathbb{C}^n$ and  we have $\ds{\frac{1}{(2\pi
t)^{n/2}} \int_{\mathbb{C}^n} |f(z)|^2 e^{-\frac{|y|^2}{2t}} dx dy
< \infty}$ $(z = x + iy).$ In this case we also have  $$ \| g\| _2
^2 = \frac{1}{(2\pi t)^{n/2}} \int_{\mathbb{C}^n} |f(z)|^2
e^{-\frac{|y|^2}{2t}} dx dy.$$

The mapping $\ds{g \rightarrow g * p_t}$ is called the
Segal-Bargmann transform and the above says that the
Segal-Bargmann transform is a unitary map from $L^2(\mathbb{R}^n)$
onto $\ds{\mathcal{O}(\mathbb{C}^n)}$ $\ds{\bigcap
L^2(\mathbb{C}^n, \mu)},$ where $\ds{d\mu(z) = \frac{1}{(2\pi
t)^{n/2}} e^{-\frac{|y|^2}{2t}} dx dy }$ and
$\mathcal{O}(\mathbb{C}^n)$ denotes the space of entire functions
on $\mathbb{C}^n.$

In the paper \cite{H}, B. C. Hall introduced a generalization of
the Segal-Bargmann transform on a compact Lie group. If $K$ is
such a group, this coherent state transform maps $L^2(K)$
isometrically onto the space of holomorphic functions in
$L^2(G,\mu_t),$ where $G$ is the complexification of $K$ and
$\mu_t$ is an appropriate heat kernel measure on $G.$ The
generalized coherent state transform is defined in terms of the
heat kernel on the compact group $K$ and its analytic continuation
to the complex group $G.$ Similar results have been proved by
various authors. See \cite{St}, \cite{HM}, \cite{HL}, \cite{KTX}
and \cite{KOS}.

Next, consider the following result on $\mathbb{R}$ due to Paley
and Wiener. A function $f \in L^2(\mathbb{R})$ admits a
holomorphic extension to the strip $\{ x + iy : |y| < t \}$ such
that $$\ds{ \sup_{|y| \leq s} \int_\mathbb{R} |f(x+iy)|^2 dx <
\infty ~ \forall s<t}$$ if and only if \bea \label{eq}
\int_{\mathbb{R}} e^{2s|\xi|} |\widetilde{f}(\xi)|^2 d\xi < \infty
~ \forall ~ s < t\eea where $\widetilde{f}$ denotes the Fourier
transform of $f.$

The condition (\ref{eq}) is the same as $$ \int_{\mathbb{R}}
|\widetilde{e^{2s{\Delta}^{\frac{1}{2}}}f}(\xi)|^2 d\xi < \infty
~~ \forall s < t$$ where $\Delta$ is the Laplacian on
$\mathbb{R}.$ This point of view was explored by R. Goodman in
Theorem 2.1 of \cite{G1}.

The condition (\ref{eq}) also equals $$ \int_{\mathbb{R}}
|e^{i(x+iy)\xi}|^2 |\widetilde{f}(\xi)|^2 d\xi < \infty ~~ \forall
y < t.$$  Here $\xi \mapsto e^{i(x+iy)\xi}$ may be seen as the
complexification of the parameters of the unitary irreducible
representations $\xi \mapsto e^{ix\xi}$ of $\mathbb{R}.$ This
point of view also was further developed by R. Goodman (see
Theorem 3.1 from \cite{G2}). Similar results were established for
the Euclidean motion group $M(2)$ of the plane $\mathbb{R}^2$ in
\cite{NS}. Aim of this paper is to prove corresponding results in
the context of general motion groups.

The plan of this paper is as follows : In the following section we
recall the representation theory and Plancherel theorem of the
motion group $M.$ We also describe the Laplacian on $M.$ In the
next section we prove the unitarity of the Segal-Bargmann
transform on $M$ and we study generalized Segal-Bargmann transform
which is an analogue of Theorem 8 and Theorem 10 in \cite{H}. The
fourth section is devoted to a study of Poisson integrals on $M$
via a Gutzmer-type formula on $M$ which is proved by using a
Gutzmer formula for compact Lie groups established by Lassalle in
1978 (see \cite{L}). This section is modelled after the work of
Goodman \cite{G1}. In the final section we prove another
characterization of functions extending holomorphically to the
complexification of $M$ which is an analogue of Theorem 3.1 of
\cite{G2}.

\section{Preliminaries}

Let $K$ be a compact, connected Lie group which acts as a linear
group on a finite dimensional real vector space $V.$ Let $M$ be
the semidirect product of $V$ and $K$ with the group law $$(x_{1},
k_{1}) \cdot (x_{2}, k_{2}) = (x_{1}+k_{1}x_{2}, k_{1}k_{2}) \txt{
where } x_1, x_2 \in V; k_1, k_2 \in K.$$

$M$ is called the motion group. Since $K$ is compact, there exists
a $K$-invariant inner product on $V$. Hence, we can assume that
$K$ is a connected subgroup of $SO(n),$ where $n = \dim V.$ When
$K=\{1\},$ $M=V\cong \mathbb{R}^{n}$ and if $K=SO(n),$ $M$ is the
Euclidean motion group. Henceforth we shall identify $V$ with
$\mathbb{R}^{n}$ and $K$ with a subgroup of $SO(n).$

The group $M$ may be identified with a matrix subgroup of
$GL(n+1,\mathbb{R})$ via the map $$ (x, k) \rightarrow \left (
\begin{matrix} k & x \\ 0 & 1 \\ \end{matrix}
\right)$$ where $x \in \mathbb{R}^{n}$ and $k \in K \subseteq
SO(n).$

We normalize the Haar measure $dm$ on $M$ such that $dm=dxdk,$
where $\ds{dx=(2\pi)^{-\frac{n}{2}}}$ $\ds{dx_1 dx_2 \cdots dx_n}$
and $dk$ is the normalized Haar measure on $K.$ Let $\mathcal{H} =
L^2(K)$ be the Hilbert space of all square integrable functions on
$K.$ Denote by $\langle \cdot , \cdot \rangle$ the Euclidean inner
product on $\mathbb{R}^{n}.$ Let $\widehat{V}$ be the dual space
of $V.$ Then we can identify $\widehat{V}$ with $\mathbb{R}^{n}$
so that $K$ acts on $\widehat{V}$ naturally by $\langle k \cdot
\xi , x \rangle = \langle \xi , k^{-1} \cdot x \rangle$ where $\xi
\in \widehat{V}, ~ x \in V, ~ k \in K.$

For any $\xi \in \widehat{V}$ let $U^{\xi}$ denote the induced
representation of $M$ by the unitary representation $x \mapsto
e^{i<\xi,x>}$ of $V.$ Then for $F \in \mathcal{H}$ and $(x,k) \in
M,$ $$ U_{(x,k)}^\xi F(u) = e^{i \langle x , u \cdot \xi \rangle}
F(k^{-1}u).$$ The representation $U^{\xi}$ is not irreducible. Any
irreducible unitary representation of $M$ is, however, contained
in $U^{\xi}$ for some $\xi \in \widehat{V}$ as an irreducible
component.

Let $K_\xi$ be the isotropy subgroup of $\xi \in \widehat{V}$ i.e.
$K_\xi = \{ k \in K : k \cdot \xi = \xi \}.$ Consider $\sigma \in
\widehat{K_{\xi}},$ the unitary dual of $K_{\xi}.$ Denote by
$\chi_{\sigma},$ $d_{\sigma}$ and $\sigma_{ij}$ the character,
degree and matrix coefficients of $\sigma$ respectively. Let $R$
be the right regular representation of $K.$ Define $$ P^{\sigma} =
d_{\sigma} \int_{K_{\xi}} \overline{\chi_{\sigma}(w)} R_w dw$$ and
$$P^{\sigma}_{\gamma} = d_{\sigma} \int_{K_{\xi}}
\overline{{\sigma_{\gamma \gamma}}(w)} R_w dw$$ where $dw$ is the
normalized Haar measure on $K_{\xi}.$ Then $P^{\sigma}$ and
$P^{\sigma}_{\gamma}$ are both orthogonal projections of
$\mathcal{H}.$ Let $\mathcal{H}^{\sigma} = P^{\sigma} \mathcal{H}$
and $\mathcal{H}^{\sigma}_{\gamma} = P^{\sigma}_{\gamma}
\mathcal{H}.$ The subspaces $\mathcal{H}^{\sigma}_{\gamma}$ are
invariant under $U^{\xi}$ for $1 \leq \gamma \leq d_{\sigma}$ and
the representations of $M$ induced on
$\mathcal{H}^{\sigma}_{\gamma}$ under $U^{\xi}$ are equivalent for
all $1 \leq \gamma \leq d_{\sigma}.$ We fix one of them and denote
it by $U^{\xi, \sigma}.$ Two representations $U^{\xi, \sigma}$ and
$U^{\xi', \sigma'}$ are equivalent if and only if there exists an
element $k \in K$ such that $\xi=k \cdot \xi'$ and $\sigma'$ is
equivalent to $\sigma^k$ where $\sigma^k(w) = \sigma(k w k^{-1})$
for $w \in K_{\xi}.$

The Mackey theory \cite{M} shows that under certain conditions on
$K$ (for details refer to Section 6.6 of \cite{F}), each $U^{\xi,
\sigma}$ is irreducible and every infinite dimensional irreducible
unitary representation is equivalent to one of $U^{\xi, \sigma}$
for some $\xi \in \mathbb{R}^n$ and $\sigma \in
\widehat{K_{\xi}}.$ Since $\ds{\mathcal{H} = \bigoplus_{\sigma \in
\widehat{K_{\xi}}} \mathcal{H}^{\sigma}}$ and
$\mathcal{H}^{\sigma} = \bigoplus_{\gamma = 1}^{d_{\sigma}}
\mathcal{H}^{\sigma}_{\gamma},$ we have $$ U^{\xi} \cong
\bigoplus_{\sigma \in \widehat{K_{\xi}}} d_{\sigma} U^{\xi,
\sigma}.$$

For any $f \in L^1(M)$ define the Fourier transform of $f$ by
$$\widehat{f}(\xi,\sigma) = \int_M f(m) U^{\xi,\sigma}_m dm.$$ Then
the Plancherel formula gives $$ \int_M |f(m)|^2 dm = \sum_{\sigma
\in \widehat{K_\xi}} d_{\sigma} \int_{\mathbb{R}^n} \|
\widehat{f}(\xi,\sigma) \|_{HS}^2 d\xi$$ where $\| \cdot \|_{HS}$
is the Hilbert-Schmidt norm of an operator. We will be working
with the generalized Fourier transform defined by
$$\widehat{f}(\xi) = \int_M f(m) U^{\xi}_m dm.$$ Then we also have  $$
\int_M |f(m)|^2 dm = \int_{\mathbb{R}^n} \| \widehat{f}(\xi)
\|_{HS}^2 d\xi.$$

Let $\underline{k}$ and $\underline{m}$ be the Lie algebras of $K$
and $M$ respectively. Then $$\underline{m} = \left \{ \left ( \begin{matrix} K & X \\ 0 & 0 \\
\end{matrix} \right) : X \in \mathbb{R}^{n}, ~ K \in  \underline{k} \right \}
.$$ Let $K_1, K_2, \cdots, K_N$ be a basis of $\underline{k}$ and
$X_1, X_2,$ $\cdots, X_n$ be a Lie algebra basis of
$\mathbb{R}^{n}.$ Define \beas
M_i &=& \left ( \begin{matrix} K_i & 0 \\ 0 & 0 \\
\end{matrix} \right) \txt { for } 1 \leq i \leq N \\
&=& \left ( \begin{matrix} 0 & X_i \\ 0 & 0 \\
\end{matrix} \right) \txt { for } N+1 \leq i \leq N+n.
\eeas Then it is easy to see that $\ds{\{ M_1, M_2, \cdots ,
M_{N+n} \}}$ forms a basis for $\underline{m}.$ The Laplacian
$\ds{\Delta_{M} = \Delta}$ is defined by $$\ds{\Delta = -(M_1^2 +
M_2^2 + \cdots + M_{N+n}^2).}$$ A simple computation using the
fact $K \subseteq SO(n)$ shows that $$\ds{\Delta = -
\Delta_{\mathbb{R}^n} - \Delta_K}$$ where
$\ds{\Delta_{\mathbb{R}^n}}$ and $\ds{\Delta_{K}}$ are the
Laplacians on $\mathbb{R}^n$ and $K$ respectively given by
$\ds{\Delta_{\mathbb{R}^n} = }$ $\ds{X_1^2 + X_2^2 + \cdots +
X_n^2}$ and $\ds{\Delta_{K} = K_1^2 + K_2^2 + \cdots + K_N^2}.$

\section{Segal-Bargmann transform and its generalisation}

Since $\ds{\Delta_{\mathbb{R}^n}}$ and $\ds{\Delta_K}$ commute, it
follows that the heat kernel $\psi_t$ associated to $\Delta$ is
given by the product of the heat kernels $p_t$ on $\mathbb{R}^n$
and $q_t$ on $K.$ In other words $$ \psi_t(x,k) = p_t(x) q_t(k) =
\frac{1}{{(4\pi t)}^{\frac{n}{2}}} e^{\frac{-|x|^2}{4t}} \sum_{\pi
\in \widehat{K}} d_\pi e^{-\frac{\lambda_{\pi}t}{2}}
\chi_{\pi}(k).$$ Here, for each unitary, irreducible
representation $\pi$ of $K,$ $d_{\pi}$ is the degree of $\pi,$
$\lambda_{\pi}$ is such that $\pi(\ds{\Delta_K}) = - \lambda_{\pi}
I$ and $\chi_{\pi}(k) = tr (\pi(k))$ is the character of $\pi.$

Denote by $G$ the complexification of $K.$ Let $\kappa_t$ be the
fundamental solution at the identity of the following equation on
$G :$ $$ \frac{du}{dt} = \frac{1}{4} \Delta_G u,$$ where
$\Delta_G$ is the Laplacian on $G.$ It should be noted that
$\kappa_t$ is the real, positive heat kernel on $G$ which is not
the same as the analytic continuation of $q_t$ on $K.$

Let $\mathcal{H}(\mathbb{C}^n \times G)$ be the Hilbert space of
holomorphic functions on $\mathbb{C}^n \times G$ which are square
integrable with respect to $\mu \bigotimes \nu (z,g)$ where
$$ d\mu(z) = \frac{1}{{(2 \pi t)}^{\frac{n}{2}}} e^{-\frac{|y|^2}{2t}} dx dy
\textmd{ on } \mathbb{C}^n $$ and $$ d\nu(g) = \int_K
\kappa_t(xg)dx \textmd{ on } G.$$ Then we have the following
theorem :

\begin{thm}
If $f \in L^2(M),$ then $f*\psi_t$ extends holomorphically to
$\mathbb{C}^n \times G.$ Moreover, the map $C_t : f \mapsto f
* \psi_t$ is a unitary map from $L^2(M)$ onto
$\mathcal{H}(\mathbb{C}^n \times G).$
\end{thm}

\begin{proof}

Let $f \in L^2(M).$ Expanding $f$ in the $K-$variable using the
Peter-Weyl theorem we obtain $$ f(x,k) = \sum_{\pi \in
\widehat{K}} d_{\pi} \sum_{i,j = 1}^{d_\pi} f_{ij}^{\pi}(x)
\phi_{ij}^{\pi}(k)$$ where for each $\pi \in \widehat{K},$ $d_\pi$
is the degree of $\pi,$ $\phi_{ij}^{\pi}$'s are the matrix
coefficients of $\pi$ and $\ds{f_{ij}^{\pi}(x) = \int_K f(x,k)
\overline{\phi_{ij}^{\pi}(k)} dk}.$ Here, the convergence is
understood in the $L^2$-sense. Moreover, by the universal property
of the complexification of a compact Lie group (see Section 3 of
\cite{H}), all the representations of $K,$ and hence all the
matrix entries, extend to $G$ holomorphically.

Since $\psi_t$ is $K$-invariant (as a function on $\mathbb{R}^n$)
a simple computation shows that $$ f * \psi_t (x, k) = \sum_{\pi
\in \widehat{K}} d_{\pi} e^{-\frac{\lambda_{\pi}t}{2}} \sum_{i,j =
1}^{d_\pi} f_{ij}^{\pi} * p_t(x) \phi_{ij}^{\pi}(k). $$ It is
easily seen that $f_{ij}^{\pi} \in L^2(\mathbb{R}^n)$ for every
$\pi \in \widehat{K}$ and $1 \leq i,j \leq d_{\pi}.$ Hence
$f_{ij}^{\pi} * p_t$ extends to a holomorphic function on
$\mathbb{C}^n$ and by the unitarity of the Segal-Bargmann
transform in $\mathbb{R}^n$ we have \bea \label{e3}
\int_{\mathbb{C}^n} |f_{ij}^{\pi} * p_t(z)|^2 \mu(y) dx dy =
\int_{\mathbb{R}^n} |f_{ij}^{\pi}(x)|^2 dx. \eea The analytic
continuation of $f
* \psi_t$ to $\mathbb{C}^n \times G$ is given by $$ f * \psi_t (z,
g) = \sum_{\pi \in \widehat{K}} d_{\pi}
e^{-\frac{\lambda_{\pi}t}{2}} \sum_{i,j = 1}^{d_\pi} f_{ij}^{\pi}
* p_t(z) \phi_{ij}^{\pi}(g).$$

We claim that the above series converges uniformly on compact
subsets of $\mathbb{C}^n \times G$ so that $f * \psi_t$ extends to
a holomorphic function on $\mathbb{C}^n \times G.$ We know from
Section 4, Proposition 1 of \cite{H} that the holomorphic
extension of the heat kernel $q_t$ on $K$ is given by
$$q_t(g) = \sum_{\pi \in \widehat{K}} d_\pi e^{-\frac{\lambda_{\pi}t}{2}}
\chi_{\pi}(g).$$ For each $g \in G,$ define the function $q_t^g(k)
= q_t(gk).$ Then $q_t^g$ is a smooth function on $K$ and is given
by \beas q_t^g(k) &=& \sum_{\pi \in \widehat{K}} d_\pi
e^{-\frac{\lambda_{\pi}t}{2}}
\chi_{\pi}(gk)\\
&=& \sum_{\pi \in \widehat{K}} d_\pi e^{-\frac{\lambda_{\pi}t}{2}}
\sum_{i,j=1}^{d_{\pi}} \phi_{ij}^{\pi}(g) \phi_{ji}^{\pi}(k).
\eeas Since $q_t^g$ is a smooth function on $K,$ we have for each
$g \in G,$ \bea \label{e1} \int_K |q_t^g(k)|^2 dk = \sum_{\pi \in
\widehat{K}} d_\pi e^{-\lambda_{\pi}t} \sum_{i,j=1}^{d_{\pi}}
|\phi_{ij}^{\pi}(g)|^2 < \infty .\eea Let $L$ be a compact set in
$\mathbb{C}^n \times G.$ For $(z,g) \in L$ we have, \bea
\label{e2} |f * \psi_t (z, g)| \leq \sum_{\pi \in \widehat{K}}
d_{\pi} e^{-\frac{\lambda_{\pi}t}{2}} \sum_{i,j = 1}^{d_\pi}
|f_{ij}^{\pi} * p_t(z)| |\phi_{ij}^{\pi}(g)|. \eea By the Fourier
inversion $$ f_{ij}^{\pi} * p_t(z) = \int_{\mathbb{R}^n}
\widetilde{f_{ij}^{\pi}} (\xi) e^{-t|\xi|^2} e^{i\xi \cdot (x +
iy)} d\xi $$ where $z = x + iy \in \mathbb{C}^n$ and
$\widetilde{f_{ij}^{\pi}}$ is the Fourier transform of
$f_{ij}^{\pi}.$ Hence, if $z$ varies in a compact subset of
$\mathbb{C}^n,$ we have \beas |f_{ij}^{\pi} * p_t(z)| &\leq&
\|f_{ij}^{\pi}\|_2 \int_{\mathbb{R}^n} e^{-2(t|\xi|^2 + y \cdot
\xi)} d\xi \\ &\leq& C \|f_{ij}^{\pi}\|_2.  \eeas Using the above
in (\ref{e2}) and applying Cauchy-Schwarz inequality we get \beas
|f * \psi_t (z, g)| &\leq& C \sum_{\pi \in \widehat{K}} d_{\pi}
\sum_{i,j = 1}^{d_\pi} \|f_{ij}^{\pi}\|_2
e^{-\frac{\lambda_{\pi}t}{2}}
|\phi_{ij}^{\pi}(g)| \\
&\leq& C \left( \sum_{\pi \in \widehat{K}} d_{\pi} \sum_{i,j =
1}^{d_\pi} \int_{\mathbb{R}^n} |f_{ij}^{\pi}(x)|^2 dx \right)
^{\frac{1}{2}} \left( \sum_{\pi \in \widehat{K}} d_{\pi} \sum_{i,j
= 1}^{d_\pi} e^{-\lambda_{\pi}t} |\phi_{ij}^{\pi}(g)|^2
\right)^{\frac{1}{2}} .\eeas Noting that $\ds{ \|f\|_2^2 =
\sum_{\pi \in \widehat{K}} d_{\pi} \sum_{i,j = 1}^{d_\pi}
\int_{\mathbb{R}^n} |f_{ij}^{\pi}(x)|^2 dx }$ and $q_t$ is a
smooth function on $G$ we prove the claim using (\ref{e1}).
Applying Theorem 2 in \cite{H} we get
$$ \int_G |f*\psi_t(z,g)|^2 d\nu(g) = \sum_{\pi \in \widehat{K}}
d_{\pi} \sum_{i,j = 1}^{d_\pi} |f_{ij}^{\pi} * p_t(z)|^2.$$
Integrating the above against $\mu(y)dx dy$ on $\mathbb{C}^n$ and
using (\ref{e3}) we obtain the isometry of $C_t$ $$
\int_{\mathbb{C}^n} \int_{G} |f
* \psi_t(z,g)|^2 \mu(y) dx dy d\nu(g) = \|f\|_2^2.$$

To prove that the map $C_t$ is surjective it suffices to prove
that the range of $C_t$ is dense in $\ds{\mathcal{H}(\mathbb{C}^n
\times G)}.$ For this, consider functions of the form $f(x,k) =
h_1(x) h_2(k) \in L^2(M)$ where $h_1 \in L^2(\mathbb{R}^n), ~~ h_2
\in L^2(K).$ Then a simple computation shows that $$ f *
\psi_t(z,g) = h_1* p_t(z) h_2*q_t(g) \txt{ for } (z,g) \in
\mathbb{C}^n \times G.$$ Suppose $\ds{F \in
\mathcal{H}(\mathbb{C}^n \times G)}$ be such that \bea \label{e4}
\int_{\mathbb{C}^n \times G} F(z,g) \overline{h_1*p_t(z)}
\overline{h_2*q_t(g)} \mu(y)dxdyd\nu(g) = 0 \eea $\forall ~  h_1
\in L^2(\mathbb{R}^n)$ and $\forall ~  h_2 \in L^2(K).$ From
(\ref{e4}) we have $$ \int_{G} \left(\int_{\mathbb{C}^n} F(z,g)
\overline{h_1*p_t(z)} d\mu(z) \right) \overline{h_2*q_t(g)}
d\nu(g) = 0,
$$ which by Theorem 2 of \cite{H} implies that $$
\int_{\mathbb{C}^n} F(z,g) \overline{h_1*p_t(z)} d\mu(z) = 0 . $$
Finally, an application of the surjectivity of Segal-Bargmann
transform on $\mathbb{R}^n$ shows that $F \equiv 0.$ Hence the
proof.
\end{proof}

In \cite{H} Brian C. Hall proved the following generalizations of
the Segal-Bargmann transfoms for $\mathbb{R}^n$ and compact Lie
groups :

\begin{thm}\label{hall}

\item [(I)] Let $\mu$ be any measurable function on $\mathbb{R}^n$
such that
\begin{itemize}
\item $\mu$ is strictly positive and locally bounded away from
zero, \item $ \ds{ \forall ~  x \in \mathbb{R}^n, ~~ \sigma(x) =
\int_{\mathbb{R}^n} e^{2x \cdot y} \mu(y) dy < \infty .}$
\end{itemize}
Define, for $z \in \mathbb{C}^n$ $$ \psi(z) = \int_{\mathbb{R}^n}
\frac{e^{ia(y)}}{\sqrt{\sigma(y)}} e^{-iy \cdot z} dy,$$ where $a$
is a real valued measurable function on $\mathbb{R}^n.$ Then the
mapping $C_{\psi} : L^2(\mathbb{R}^n) \rightarrow
\mathcal{O}(\mathbb{C}^n)$ defined by $$ C_{\psi}(z) =
\int_{\mathbb{R}^n} f(x) \psi(z-x) dx $$ is an isometric
isomorphism of $L^2(\mathbb{R}^n)$ onto $\mathcal{O}(\mathbb{C}^n)
\bigcap L^2(\mathbb{C}^n,dx \mu(y)dy).$

\item [(II)] Let $K$ be a compact Lie group and $G$ be its
complexification. Let $\nu$ be a measure on $G$ such that
\begin{itemize}
\item $\nu$ is bi-$K$-invariant, \item $\nu$ is given by a
positive density which is locally bounded away from zero, \item
For each irreducible representation $\pi$ of $K,$ analytically
continued to $G,$ $$ \delta(\pi) = \frac{1}{dim V_{\pi}} \int_{G}
\|\pi(g^{-1})\|^2 d\nu(g) < \infty.$$
\end{itemize}
Define $ \ds {\tau (g) = \sum _ {\pi \in \widehat{K}}
\frac{d_{\pi}}{\sqrt{\delta(\pi)}} Tr (\pi(g^{-1}) U_{\pi})} $
where $g \in G$ and $U_{\pi}$'s are arbitrary unitary matrices.
Then the mapping $$C_\tau f(g) = \int_{K} f(k) \tau(k^{-1}g)dk$$
is an isometric isomorphism of $L^2(K)$ onto $\ds{\mathcal{O}(G)
\bigcap L^2(G, d\nu(w)).}$

\end{thm}

A similar result holds for $M$. Let $\mu$ be any real-valued
$K$-invariant function on $\mathbb{R}^n$ such that it satisfies
the conditions of Theorem \ref{hall} (I). Define, for $z \in
\mathbb{C}^n$ $$ \psi(z) = \int_{\mathbb{R}^n}
\frac{e^{ia(y)}}{\sqrt{\sigma(y)}} e^{-iy.z} dy,$$ where $a$ is a
real valued measurable K-invariant function on $\mathbb{R}^n.$
Next, let $\nu, ~ \delta $ and $\tau$ be as in Theorem \ref{hall}
(II). Also define $\ds{ \phi(z,g) = \psi(z) \tau(g) }$ for $z \in
\mathbb{C}^n, ~ g \in G.$ It is easy to see that $\phi(z,w)$ is a
holomorphic function on $\mathbb{C}^n \times G.$ Then it is easy
to prove the following analogue of Theorem \ref{hall} for $M.$

\begin{thm}
The mapping $$\ds{C_\phi f(z,g) = \int_{M} f(\xi, k)
\phi((\xi,k)^{-1}(z,g))d\xi dk}$$ is an isometric isomorphism of
$L^2(M)$ onto $$\ds{\mathcal{O}(\mathbb{C}^n \times G) \bigcap
L^2(\mathbb{C}^n \times G, \mu(y)dxdyd\nu(g)).}$$
\end{thm}

\section{Gutzmer's formula and Poisson Integrals}

In this section first we briefly recall Gutzmer's formula on
compact, connected Lie groups given by Lassalle in \cite{L}. Let
$\underline{k}$ and $\underline{g}$ be the Lie algebras of a
compact, connected Lie group $K$ and its complexification $G.$
Then we can write $\underline{g} = \underline{k} + \underline{p}$
where $\underline{p}=i\underline{k}$ and any element $g \in G$ can
be written in the form $g=k\exp{iH}$ for some $k \in K, ~ H \in
\underline{k}.$ If $\underline{h}$ is a maximal, abelian
subalgebra of $\underline{k}$ and $\underline{a}=i\underline{h}$
then every element of $\underline{p}$ is conjugate under $K$ to an
element of $\underline{a}.$ Thus each $g \in G$ can be written
(non-uniquely) in the form $g = k_1 \exp {(iH)} k_2 $ for $k_1,
k_2 \in K$ and $H \in \underline{h}.$ If $ k_1 \exp {(iH_1)} k_1'
= k_2 \exp {(iH_2)} k_2', $ then there exists $w \in W,$ the Weyl
group with respect to $\underline{h},$ such that $H_1 = w \cdot
H_2$ where $\cdot$ denotes the action of the Weyl group on
$\underline{h}.$ Since $K$ is compact, there exists an Ad-$
K$-invariant inner product on $\underline{k},$ and hence on
$\underline{g}.$ Let $|\cdot|$ denote the norm with respect to the
said inner product. Then we have the following Gutzmer's formula
by Lassalle.

\begin{thm}\label{lassalle}
Let $f$ be holomorphic in $K\exp(i\Omega_r)K \subseteq G$ where
$\Omega_r = \{ H \in \underline{k} : |H| <r \}.$ Then we have $$
\int_K \int_K |f(k_1 \exp{iH} k_2)|^2 dk_1 dk_2 = \sum_{\pi \in
\widehat{K}} \|\widehat{f}(\pi)\|_{HS}^2 \chi_{\pi}(\exp{2iH})$$
where $H \in \Omega_r$ and $\widehat{f}(\pi)$ is the
operator-valued Fourier transform of $f$ at $\pi$ defined by
$\ds{\widehat{f}(\pi) = \int_K f(k) \pi(k^{-1}) dk.}$
\end{thm}

For the proof of above see \cite{L}. We prove a Gutzmer-type
result on $M$ using Lassalle's theorem above. Define
$\ds{\Omega_{t,r} = \{ (z,g) \in }$ $\ds{ \mathbb{C}^n \times G :
|Im z|<t, |H|<r \txt{ where } g = }$ $\ds{k_1 \exp{iH} k_2, k_1,
k_2}$ $\ds{\in K, H \in \underline{h} \}}.$ Notice that the domain
$\Omega_{t,r}$ is well defined since $|\cdot|$ is invariant under
the Weyl group action.

\begin{lem}\label{lemma}
Let $f \in L^2(M)$ extend holomorphically to the domain
$\Omega_{t,r}$ and $$ \sup_{\left\{ |y|<s, ~ |H|<q \right\}}
\int_{\mathbb{R}^n} \int_K \int_K |f(x+iy,k_1\exp{(iH)}k_2)|^2
dk_1 dk_2 dx < \infty $$ $\forall ~ s<t$ and $q<r.$ Then \beas
\int_{\mathbb{R}^n} \int_K \int_K |f(x+iy,k_1\exp{(iH)}k_2)|^2
dk_1 dk_2 dx \\ = \sum_{\pi \in \widehat{K}} d_{\pi} \sum_{i,j =
1}^{d_\pi} \left(\int_{\mathbb{R}^n}
|\widetilde{f_{ij}^{\pi}}(\xi)|^2 e ^{-2 \xi.y} d\xi \right)
\chi_{\pi}(\exp{2iH})\eeas provided $|y|<t$ and $|H|<r.$
Conversely, if $$ \sup_{\left\{ |y|<s, ~ |H|<q \right\}} \sum_{\pi
\in \widehat{K}} d_{\pi} \sum_{i,j = 1}^{d_\pi}
\left(\int_{\mathbb{R}^n} |\widetilde{f_{ij}^{\pi}}(\xi)|^2 e ^{-2
\xi.y} d\xi \right) \chi_{\pi}(\exp{2iH}) < \infty ~ \forall ~ s<t
\txt{ and } q<r $$ then $f$ extends holomorphically to the domain
$\ds{\Omega_{t,r}}$ and $$ \sup_{\left\{ |y|<s, ~ |H|<q \right\}}
\int_{\mathbb{R}^n} \int_K \int_K |f(x+iy,k_1\exp{(iH)}k_2)|^2
dk_1 dk_2 dx < \infty $$ $\forall ~ s<t$ and $q<r.$
\end{lem}

\begin{proof}
Notice that $\ds{f_{ij}^{\pi}(x) = \int_K f(x,k)
\overline{\phi_{ij}^{\pi}(k)} dk}.$ It follows that $f_{ij}^{\pi}$
has a holomorphic extension to $\{ z \in \mathbb{C}^n : |Im z| < t
\}$ and $$\ds{ \sup_{|y|<s} \int_{\mathbb{R}^n}
|f_{ij}^{\pi}(x+iy)|^2 dx < \infty ~ \forall ~  s<t. }$$
Consequently, $$ \int_{\mathbb{R}^n} |f_{ij}^{\pi}(x+iy)|^2 dx =
\int_{\mathbb{R}^n} |\widetilde{f_{ij}^{\pi}}(\xi)|^2 e ^{-2
\xi.y} d\xi \txt{ for } |y|<s ~ \forall ~  s<t. $$ Now, for each
fixed $z \in \mathbb{C}^n$ with $|Im z|<s$ the function $g \ra
f(z,g)$ is holomorphic in the domain $\ds{\{ g \in G : |H|<r \txt{
where } g = k_1 \exp{iH} k_2, k_1, k_2 \in K, H \in \underline{h}
\}}$ for every $s<t$ and $q<r$ and so admits a holomorphic Fourier
series (as in \cite{H}) $$ f(z,g) = \sum_{\pi \in \widehat{K}}
d_{\pi} \sum_{i,j = 1}^{d_\pi} a_{ij}^{\pi}(z)
\phi_{ij}^{\pi}(g).$$ It follows that $a_{ij}^{\pi}(z) =
f_{ij}^{\pi}(z)$ for every $\pi \in \widehat{K}$ and $1 \leq i,j
\leq d_{\pi}.$ Hence by using Theorem \ref{lassalle} we have for
$(z,g) \in \Omega_{t,r},$ \beas \int_K \int_K
|f(x+iy,k_1\exp{(iH)}k_2)|^2 dk_1 dk_2 &=& \sum_{\pi \in
\widehat{K}} \| \widehat{f_z}(\pi) \|_{HS}^2 \chi_{\pi}(\exp{2iH})\\
&=& \sum_{\pi \in \widehat{K}} \sum_{i,j = 1}^{d_\pi}
|f_{ij}^{\pi}(z)|^2 \chi_{\pi}(\exp{2iH}) \eeas where $f_z(g) =
f(z,g).$ Integrating over $\mathbb{R}^n$ we get \beas &&
\int_{\mathbb{R}^n} \int_K \int_K |f(x+iy,k_1\exp{(iH)}k_2)|^2
dk_1 dk_2 dx \\ &=& \sum_{\pi \in \widehat{K}} \sum_{i,j =
1}^{d_\pi} \int_{\mathbb{R}^n} |f_{ij}^{\pi}(x + iy)|^2 dx ~~
\chi_{\pi}(\exp{2iH}) \\ &=& \sum_{\pi \in \widehat{K}} d_{\pi}
\sum_{i,j = 1}^{d_\pi} \left(\int_{\mathbb{R}^n}
|\widetilde{f_{ij}^{\pi}}(\xi)|^2 e ^{-2 \xi.y} d\xi \right)
\chi_{\pi}(\exp{2iH}).\eeas Hence the first part of the lemma is
proved. Converse can also be proved similarly.

\end{proof}

Recall that the Laplacian $\Delta$ on $M$ is given by $\ds{\Delta
= -\Delta_{\mathbb{R}^n} - \Delta_K.}$ If $f \in L^2(M)$ it is
easy to see that $$ e^{-t\Delta^{\frac{1}{2}}} f (x,k) = \sum_{\pi
\in \widehat{K}} d_{\pi} \sum_{i,j = 1}^{d_\pi}
\left(\int_{\mathbb{R}^n} e^{-t(|\xi|^2 +
\lambda_{\pi})^{\frac{1}{2}}} \widetilde{f_{ij}^{\pi}}(\xi) e
^{i\xi \cdot y} d\xi \right) \phi_{ij}^{\pi}(k).$$ We have the
following (almost) characterization of the Poisson integrals. Let
$\Omega_{t,r}$ denote the domain defined in Lemma \ref{lemma}.

\begin{thm}\label{theorem}

Let $f \in L^2(M).$ Then there exists a constant $N$ such that $g
= e^{-t\Delta^{\frac{1}{2}}} f$ extends to a holomorphic function
on the domain $\Omega_{\frac{t}{\sqrt{2}},\frac{t\sqrt{2}}{N}}$
and $$ \sup_{\left\{|y|<\frac{t}{\sqrt{2}}, ~ |H| \leq
\frac{t\sqrt{2}}{N} \right\}} \int_{\mathbb{R}^n} \int_K \int_K
|g(x+iy,k_1\exp{(iH)}k_2)|^2 dk_1 dk_2 dx < \infty .$$ Conversely,
there exists a fixed constant $C$ such that whenever $g$ is a
holomorphic function on $\Omega_{t,\frac{2t}{C}}$ and
$$ \ds{\sup_{\left\{|y|<s, ~ |H|< \frac{2s}{C} \right\}} \int_{\mathbb{R}^n} \int_K
\int_K |g(x+iy,k_1\exp{(iH)}k_2)|^2 dk_1 dk_2 dx < \infty \txt{
for } s<t,}$$ then $\forall ~  s<t, ~ \exists f \in L^2(M)$ such
that $e^{-s\Delta^{\frac{1}{2}}} f = g.$

\end{thm}

\begin{proof}
We know that, if $f \in L^2(M)$ then $$ g (x,k) =
e^{-t\Delta^{\frac{1}{2}}} f (x,k) = \sum_{\pi \in \widehat{K}}
d_{\pi} \sum_{i,j = 1}^{d_\pi} \left(\int_{\mathbb{R}^n}
e^{-t(|\xi|^2 + \lambda_{\pi})^{\frac{1}{2}}}
\widetilde{f_{ij}^{\pi}}(\xi) e ^{i\xi \cdot y} d\xi \right)
\phi_{ij}^{\pi}(k).$$ Also, $\ds{g(x,k) =\sum_{\pi \in
\widehat{K}} d_{\pi} \sum_{i,j = 1}^{d_\pi} g_{ij}^{\pi}(x)
\phi_{ij}^{\pi}(k) }$ with $\ds{ \widetilde{g_{ij}^{\pi}}(\xi) =
\widetilde{f_{ij}^{\pi}}(\xi) e^{-t(|\xi|^2 +
\lambda_{\pi})^{\frac{1}{2}}}}.$ From Lemma 6 and 7 of \cite{H} we
know that there exist constants $a$, $M$ such that $\lambda_{\pi}
\geq a|\mu|^2$ and $|\chi_{\pi}(\exp{iY})| \leq d_{\pi}
e^{M|Y||\mu|}$ where $\mu$ is the highest weight of $\pi.$ Hence
we have
$$ |\chi_{\pi}(\exp 2iH)| \leq d_{\pi} e^{2M|H||\mu|} \leq d_{\pi}
e^{N|H| \sqrt{\lambda_{\pi}}}$$ where $N=\frac{2M}{\sqrt{a}}.$ If
$s \leq \frac{t\sqrt{2}}{N}$ it is easy to show that $$
\sup_{\left\{\xi \in \mathbb{R}^n, ~ \lambda_{\pi} \geq 0
\right\}} e^{-2t(|\xi|^2 + \lambda_{\pi})^{\frac{1}{2}}} e^{2|\xi|
|y|} e^{N|\sqrt{\lambda_{\pi}}|s} \leq C < \infty \txt{ for } |y|
\leq \frac{t}{\sqrt{2}}.
$$ It follows that $$ \sup_{\left\{|y|<\frac{t}{\sqrt{2}}, ~ |H|
\leq \frac{t\sqrt{2}}{N} \right\}} \sum_{\pi \in \widehat{K}}
d_{\pi} \sum_{i,j = 1}^{d_\pi} \left(\int_{\mathbb{R}^n}
|\widetilde{g_{ij}^{\pi}}(\xi)|^2 e ^{-2\xi \cdot y} d\xi \right)
e^{N\sqrt{\lambda_{\pi}}|H|} < \infty.$$ So we have $$
\sup_{\left\{|y|<\frac{t}{\sqrt{2}}, ~ |H| \leq
\frac{t\sqrt{2}}{N} \right\}} \sum_{\pi \in \widehat{K}} d_{\pi}
\sum_{i,j = 1}^{d_\pi} \left(\int_{\mathbb{R}^n}
|\widetilde{g_{ij}^{\pi}}(\xi)|^2 e ^{-2\xi \cdot y} d\xi \right)
\chi_{\pi}(\exp 2iH) < \infty.$$ Hence by Lemma \ref{lemma} we
have proved one part of the theorem.

To prove the converse, we first show that there exist constants
$A,C$ such that \bea \label{e7} \int_{|H|=r} \chi_{\pi} (\exp
{2iH}) d\sigma_r(H) \geq d_{\pi} A e^{Cr\sqrt{\lambda_{\pi}}}\eea
where $d\sigma_r(H)$ is the normalized surface measure on the
sphere $\{ H \in \underline{h} : |H| = r\} \subseteq \mathbb{R}^m$
where $m = \dim \underline{h}.$ If $H \in \underline{a},$ then
there exists a non-singular matrix $Q$ and pure-imaginary valued
linear forms $\nu_1, \nu_2, \cdots , \nu_{d_{\pi}}$ on
$\underline{a}$ such that $$ Q \pi(H) Q^{-1} = diag(\nu_1(H),
\nu_2(H), \cdots , \nu_{d_{\pi}}(H))$$ where $diag(a_1, a_2,
\cdots , a_k)$ denotes $k \times k$ order diagonal matrix with
diagonal entries $a_1, a_2, \cdots , a_k.$ Now, $\nu(H)=i\langle
\nu, H\rangle$ where $\nu$ is a weight of $\pi.$ Then $$ \exp (2i
Q \pi(H) Q^{-1}) = Q \exp(2i\pi(H)) Q^{-1} = diag(e^{2i\nu_1(H)},
e^{2i\nu_2(H)}, \cdots , e^{2i\nu_{d_{\pi}}(H)}).$$ Hence \beas
\chi_{\pi}(\exp2iH) &=& Tr (Q \exp(2i\pi(H)) Q^{-1}) \\ &=&
e^{-2\langle \nu_1, H\rangle} + e^{-2\langle \nu_2, H\rangle} +
\cdots + e^{-2\langle \nu_{d_{\pi}}, H\rangle} \\ &\geq&
e^{-2\langle \mu, H\rangle}\eeas where $\mu$ is the highest weight
corresponding to $\pi.$ Integrating the above over $|H|=r$ we get
\beas \int_{|H|=r} \chi_{\pi} (\exp {2iH}) d\sigma_r(H) &\geq&
\int_{|H|=r}  e^{-2\langle \mu, H\rangle} d\sigma_r(H)\\
&=& \frac{J_{\frac{m}{2}-1}(2ir|\mu|)}{(2ir|\mu|)^{\frac{m}{2}-1}}
\\ &\geq& B e^{r|\mu|}\eeas where $J_{\frac{m}{2}-1}$ is the Bessel function of
order ${\frac{m}{2}-1}.$ By Weyl's dimension formula we know that
$d_{\pi}$ can be written as a polynomial in $\mu$ and
$\lambda_{\pi} \approx |\mu|^2.$ Hence we have
$$ \int_{|H|=r} \chi_{\pi} (\exp {2iH}) d\sigma_r(H) \geq A ~
d_{\pi} e^{Cr\sqrt{\lambda_{\pi}}}$$ for some $C.$ Consider the
domain $\Omega_{t,\frac{2t}{C}}$ for this $C.$ Let $g$ be a
holomorphic function on $\Omega_{t,\frac{2t}{C}}$ and
$$ \ds{\sup_{\left\{|y|<s, ~ |H|< \frac{2s}{C} \right\}}
\int_{\mathbb{R}^n} \int_K \int_K |g(x+iy,k_1\exp{(iH)}k_2)|^2
dk_1 dk_2 dx < \infty \txt{ for } s<t.}$$ By Lemma \ref{lemma} we
have $$\ds{ \sup_{\left\{|y|<s, ~ |H|< \frac{2s}{C} \right\}}
\sum_{\pi \in \widehat{K}} d_{\pi} \sum_{i,j = 1}^{d_\pi}
\left(\int_{\mathbb{R}^n} |\widetilde{g_{ij}^{\pi}}(\xi)|^2 e ^{-2
\xi.y} d\xi \right) \chi_{\pi}(\exp{2iH}) < \infty ~ \forall ~
s<t.}$$ Integrating the above over $|H|=r = \frac{2s}{C}$ and
$|y|=s<t$ we obtain $$ \sum_{\pi \in \widehat{K}} d_{\pi}
\sum_{i,j = 1}^{d_\pi} \left(\int_{\mathbb{R}^n}
|\widetilde{g_{ij}^{\pi}}(\xi)|^2
\frac{J_{\frac{n}{2}-1}(2is|\xi|)}{(2is|\xi|)^{\frac{n}{2}-1}}
d\xi \right) \int_{|H|=r} \chi_{\pi}(\exp{2iH}) d\sigma_r(H) <
\infty.$$ Noting that
$\ds{\frac{J_{\frac{n}{2}-1}(2is|\xi|)}{(2is|\xi|)^{\frac{n}{2}-1}}
\sim e^{2s|\xi|}}$ for large $|\xi|$ and using (\ref{e7}) we
obtain
$$ \sum_{\pi \in \widehat{K}} d_{\pi} \sum_{i,j = 1}^{d_\pi}
\int_{\mathbb{R}^n} |\widetilde{g_{ij}^{\pi}}(\xi)|^2 e^{2s|\xi|}
e^{2s\sqrt{\lambda_{\pi}}} d\xi < \infty \txt{ for } s<t. $$ This
surely implies that $$ \sum_{\pi \in \widehat{K}} d_{\pi}
\sum_{i,j = 1}^{d_\pi} \int_{\mathbb{R}^n}
|\widetilde{g_{ij}^{\pi}}(\xi)|^2 e^{2s(|\xi|^2 +
\lambda_{\pi})^{\frac{1}{2}}} d\xi < \infty \txt{ for } s<t.
$$ Defining $\widetilde{f_{ij}^{\pi}}(\xi)$ by $\widetilde{f_{ij}^{\pi}}(\xi) =
\widetilde{g_{ij}^{\pi}}(\xi) e^{s(|\xi|^2 +
\lambda_{\pi})^{\frac{1}{2}}}$ we obtain
$$f(x,k) = \sum_{\pi \in \widehat{K}} d_{\pi} \sum_{i,j =
1}^{d_\pi} f_{ij}^{\pi}(x) \phi_{ij}^{\pi}(k) \in L^2(M)$$ and $g
= e^{-s\Delta^{\frac{1}{2}}} f.$

\end{proof}

\section{Complexified representations and Paley-Wiener type theorems}

Recall the representations $U^{\xi}$ and the generalized Fourier
transform $\widehat{f}(\xi)$ from the introduction where
$$\widehat{f}(\xi) = \int_M f(m) U^{\xi}_m dm.$$ For $(x,k) \in M$
and matrix coefficients $\phi_{ij}^{\pi}$ of $\pi$ we have $$
\left( U_{(x,k)}^{\xi} \phi_{ij}^{\pi} \right) (u) = e^{i\langle
x,u \cdot \xi \rangle} \phi_{ij}^{\pi}(k^{-1}u).$$ This action of
$U_{(x,k)}^{\xi}$ on $\phi_{ij}^{\pi}$ can clearly be analytically
continued to $\mathbb{C}^n \times G$ and we obtain
$$ \left( U_{(z,g)}^{\xi} \phi_{ij}^{\pi} \right) (u) = e^{i\langle
x,u \cdot \xi \rangle} e^{-\langle y,u \cdot \xi \rangle}
\phi_{ij}^{\pi}(e^{-iH} k^{-1} u)
$$ where $(z,g) \in \mathbb{C}^n \times G$ and $z=x+iy
\in \mathbb{C}^n$ and $g = ke^{iH} \in G.$

We also note that the action of $K \subseteq SO(n)$ on
$\mathbb{R}^n$ naturally extends to an action of $G \subseteq
SO(n,\mathbb{C})$ on $\mathbb{C}^n.$ Then we have the following
theorem:

\begin{thm}\label{thm}
Let $f \in L^2(M).$ Then $f$ extends holomorphically to
$\mathbb{C}^n \times G$ with $$\int_{|y| = r} \int_K
\int_{\mathbb{R}^n} |f(e^{-iH}(x+iy),e^{-iH}k)|^2 dx dk d\mu_r(y)
< \infty ~~~~~~ \forall ~~ H \in \underline{k}$$ (where $\mu_r$ is
the normalized surface area measure on the sphere $\{ |y| = r \}
\subseteq \mathbb{R}^n$) iff
$$\ds{\int_{\mathbb{R}^n} \int_{|y| = r} \|U_{(z,g)}^{\xi}
\widehat{f}(\xi) \|_{HS}^2 d\mu_r(y) d\xi} < \infty $$ where
$z=x+iy \in \mathbb{C}^n,$ $g = ke^{iH} \in G.$ In this case we
also have \beas && \int_{\mathbb{R}^n} \int_{|y| = r}
\|U_{(z,g)}^{\xi} \widehat{f}(\xi) \|_{HS}^2 d\mu_r(y) d\xi \\ &=&
\int_{|y| = r} \int_K \int_{\mathbb{R}^n}
|f(e^{-iH}(x+iy),e^{-iH}k)|^2 dx dk d\mu_r(y).\eeas
\end{thm}

\vspace*{0.3in}

We know that any $f \in L^2(M)$ can be expanded in the $K$
variable using the Peter Weyl theorem to obtain \bea \label{eqn3}
f(x,k) = \sum_{\pi \in \widehat{K}} d_{\pi} \sum_{i,j = 1}^{d_\pi}
f_{ij}^{\pi}(x) \overline{\phi_{ij}^{\pi}(k)}\eea where for each
$\pi \in \widehat{K},$ $d_\pi$ is the degree of $\pi,$
$\phi_{ij}^{\pi}$'s are the matrix coefficients of $\pi$ and
$\ds{f_{ij}^{\pi}(x) = \int_K f(x,k) \phi_{ij}^{\pi}(k) dk}.$

Now, for $F \in L^2(\mathbb{R}^n),$ consider the decomposition of
the function $k \mapsto F(k \cdot x)$ in terms of the irreducible
unitary representations of $K$ given by
$$F(k \cdot x) = \sum_{\lambda \in \widehat{K}} d_{\lambda} \sum_{l,m=1}^{d_{\lambda}}
F^{lm}_{\lambda}(x) \phi_{lm}^{\lambda}(k)$$ where
$\ds{F^{lm}_{\lambda}(x) = \int_K F(k \cdot x)
\overline{\phi_{lm}^{\lambda} (k)} dk}.$ Putting $k = e,$ the
identity element of $K,$ we obtain
$$F(x) = \sum_{\lambda \in \widehat{K}} d_{\lambda}
\sum_{l=1}^{d_{\lambda}} F^{ll}_{\lambda}(x).$$ Then it is easy to
see that for $u \in K,$ \bea \label{eqn1} F^{ll}_{\lambda}(u \cdot
x) = \sum_{m=1}^{d_{\lambda}} F^{lm}_{\lambda}(x)
\phi_{lm}^{\lambda}(u).\eea It also follows that the Euclidean
Fourier transform $\widetilde{F^{lm}_{\lambda}}$ of
$F^{lm}_{\lambda}$ satisfies \bea \label{eqn2}
\widetilde{F^{ll}_{\lambda}}(u \cdot x) = \sum_{m=1}^{d_{\lambda}}
\phi_{lm}^{\lambda}(u) \widetilde{F^{lm}_{\lambda}}(x) ~~~~~
\forall ~~ u \in K.\eea From the above and the fact that
$f_{ij}^{\pi} \in L^2(\mathbb{R}^n)$ for every $\pi \in
\widehat{K}$ and $1 \leq i,j \leq d_{\pi}$ it follows that any $f
\in L^2(M)$ can be written as
$$ f(x,k) = \sum_{\pi \in \widehat{K}} d_{\pi} \sum_{\lambda \in
\widehat{K}} d_{\lambda} \sum_{i,j=1}^{d_{\pi}}
\sum_{l=1}^{d_{\lambda}} (f_{ij}^{\pi})^{ll}_{\lambda}(x)
\overline{\phi_{ij}^{\pi}(k)}.$$ We need the following lemma to
prove Theorem \ref{thm}:

\begin{lem}\label{lem}
For fixed $\pi, \lambda \in \widehat{K},$ the theorem is true for
functions of the form $$f(x,k) = \sum_{i,j=1}^{d_{\pi}}
\sum_{l=1}^{d_{\lambda}} f^{ll}_{ij}(x)
\overline{\phi_{ij}^{\pi}(k)}$$ where for simplicity we write
$(f_{ij}^{\pi})^{ll}_{\lambda}$ as $f^{ll}_{ij}.$
\end{lem}

\begin{proof}

For $\xi \in \mathbb{R}^n,$ $u \in K,$ $\gamma \in \widehat{K}$
and $1 \leq p,q \leq d_{\gamma}$ we have \beas \left(
\widehat{f}(\xi) \overline{\phi_{pq}^{\gamma}} \right) (u) &=&
\int_{\mathbb{R}^n} \int_K \sum_{i,j=1}^{d_{\pi}}
\sum_{l=1}^{d_{\lambda}} f^{ll}_{ij}(x)
\overline{\phi_{ij}^{\pi}(k)} e^{i \langle x , u \cdot \xi
\rangle} \overline{\phi_{pq}^{\gamma}(k^{-1}u)} dk dx \\
&=& \sum_{i,j=1}^{d_{\pi}} \sum_{l=1}^{d_{\lambda}}
\widetilde{f^{ll}_{ij}}(u \cdot \xi) \sum_{t=1}^{d_{\gamma}}
\phi_{qt}^{\gamma}(u^{-1}) \langle \phi_{ij}^{\pi} ,
\phi_{tp}^{\gamma} \rangle_{L^2(K)} \\
&=& \frac{\delta_{\gamma \pi}}{d_{\pi}} \sum_{i=1}^{d_{\pi}}
\sum_{l,m=1}^{d_{\lambda}} \widetilde{f^{lm}_{ip}}(\xi)
\phi_{lm}^{\lambda}(u) \phi_{qi}^{\pi}(u^{-1}) \eeas by
(\ref{eqn2}) and Schur's orthogonality relations where
$\delta_{\gamma \pi}$ is the Kronecker delta in the sense of
equivalence of unitary representations. Then we have $$ \left(
U_{(x+iy,ke^{iH})}^{\xi} \widehat{f}(\xi)
\overline{\phi_{pq}^{\gamma}} \right)(u) = \frac{\delta_{\gamma
\pi}}{d_{\pi}} e^{i \langle x+iy, u \cdot \xi \rangle}
\sum_{i=1}^{d_{\pi}} \sum_{l,m=1}^{d_{\lambda}}
\widetilde{f^{lm}_{ip}}(\xi) \phi_{lm}^{\lambda}(e^{-iH}k^{-1}u)
\phi_{qi}^{\pi}(u^{-1}ke^{iH}).$$ Hence \beas && \|
U_{(x+iy,ke^{iH})}^{\xi} \widehat{f}(\xi) \|_{HS}^2 \\
&=& \frac{1}{d_{\pi}} \sum_{p,q=1}^{d_{\pi}} \int_K e^{-2 \langle
y, u \cdot \xi \rangle} \left| \sum_{i=1}^{d_{\pi}}
\sum_{l,m=1}^{d_{\lambda}} \widetilde{f^{lm}_{ip}}(\xi)
\phi_{lm}^{\lambda}(e^{-iH}k^{-1}u) \phi_{qi}^{\pi}(u^{-1}ke^{iH})
\right|^2 du. \eeas Integrating the above over $|y|=r,$ we obtain
\bea \label{eqn4} \left. \begin{array}{rcll} && \ds{ \int_{|y|=r}
\| U_{(x+iy,ke^{iH})}^{\xi} \widehat{f}(\xi) \|_{HS}^2 d\mu_r(y)} \\
&=& \ds{ \frac{1}{d_{\pi}}
\frac{J_{\frac{n}{2}-1}(2ir|\xi|)}{(2ir|\xi|)^{\frac{n}{2}-1}}
\sum_{p,q=1}^{d_{\pi}} \int_K \left| \sum_{i=1}^{d_{\pi}}
\sum_{l,m=1}^{d_{\lambda}} \widetilde{f^{lm}_{ip}}(\xi)
\phi_{lm}^{\lambda}(e^{-iH}u) \phi_{qi}^{\pi}(u^{-1}e^{iH})
\right|^2 du} \end{array}\right. \eea where $J_{\frac{n}{2}-1}$ is
the Bessel function of order $\frac{n}{2}-1$ and $\mu_r$ is the
normalized surface area measure on the sphere $\{ |y| = r \}
\subseteq \mathbb{R}^n.$

Let $\mathcal{H}_{\pi}$ be the Hilbert space on which $\pi(k)$
acts for every $k \in K$ and $e_1, e_2, \cdots , e_{d_{\pi}}$ be a
basis of $\mathcal{H}_{\pi}.$ Then, for any $c_i, 1 \leq i \leq
d_{\pi},$ \beas && \sum_{q=1}^{d_{\pi}} \left|
\sum_{i=1}^{d_{\pi}} c_i \phi_{qi}^{\pi}(u^{-1}e^{iH}) \right|^2 \\
&=& \sum_{q=1}^{d_{\pi}} \sum_{i=1}^{d_{\pi}} c_i
\phi_{qi}^{\pi}(u^{-1}e^{iH}) \sum_{a=1}^{d_{\pi}} \overline{c_a}
\overline{\phi_{qa}^{\pi}(u^{-1}e^{iH})} \\
&=& \sum_{i,a=1}^{d_{\pi}} c_i \overline{c_a} \sum_{q=1}^{d_{\pi}}
\left\langle \pi(u^{-1}e^{iH})e_i ,e_q \right\rangle \left\langle
e_q , \pi(u^{-1}e^{iH})e_a \right\rangle \\
&=& \sum_{i,a=1}^{d_{\pi}} c_i \overline{c_a} \left\langle
\pi(u^{-1})\pi(e^{iH})e_i , \pi(u^{-1})\pi(e^{iH})e_a \right\rangle \\
&=& \sum_{q=1}^{d_{\pi}} \left| \sum_{i=1}^{d_{\pi}} c_i
\phi_{qi}^{\pi}(e^{iH}) \right|^2, \eeas since $\pi$ is a unitary
representation of $K.$ So, we have \beas \sum_{q=1}^{d_{\pi}}
\left| \sum_{i=1}^{d_{\pi}} \sum_{l,m=1}^{d_{\lambda}}
\widetilde{f^{lm}_{ip}}(\xi) \phi_{lm}^{\lambda}(e^{-iH}u)
\phi_{qi}^{\pi}(u^{-1}e^{iH}) \right|^2 = \sum_{q=1}^{d_{\pi}}
\left| \sum_{i=1}^{d_{\pi}} \sum_{l,m=1}^{d_{\lambda}}
\widetilde{f^{lm}_{ip}}(\xi) \phi_{lm}^{\lambda}(e^{-iH}u)
\phi_{qi}^{\pi}(e^{iH}) \right|^2. \eeas  Hence from (\ref{eqn4})
we get that \beas && \int_{|y|=r}
\| U_{(x+iy,ke^{iH})}^{\xi} \widehat{f}(\xi) \|_{HS}^2 d\mu_r(y) \\
&=& \ds{ \frac{1}{d_{\pi}}
\frac{J_{\frac{n}{2}-1}(2ir|\xi|)}{(2ir|\xi|)^{\frac{n}{2}-1}}
\sum_{p,q=1}^{d_{\pi}} \int_K \left| \sum_{i=1}^{d_{\pi}}
\sum_{l,m,k=1}^{d_{\lambda}} \widetilde{f^{lm}_{ip}}(\xi)
\phi_{lk}^{\lambda}(e^{-iH}) \phi_{km}^{\lambda}(u)
\phi_{qi}^{\pi}(e^{iH}) \right|^2 du} \\
&=& \frac{1}{d_{\pi} d_{\lambda}}
\frac{J_{\frac{n}{2}-1}(2ir|\xi|)}{(2ir|\xi|)^{\frac{n}{2}-1}}
\sum_{p,q=1}^{d_{\pi}} \sum_{m,k=1}^{d_{\lambda}} \left|
\sum_{i=1}^{d_{\pi}} \sum_{l=1}^{d_{\lambda}}
\widetilde{f^{lm}_{ip}}(\xi) \phi_{lk}^{\lambda}(e^{-iH})
\phi_{qi}^{\pi}(e^{iH}) \right|^2, \eeas by Schur's orthogonality
relations. The above can also be written as \bea \label{eqn5}
\left. \begin{array}{rcll} && \ds{ \int_{|y|=r}
\| U_{(x+iy,ke^{iH})}^{\xi} \widehat{f}(\xi) \|_{HS}^2 d\mu_r(y)} \\
&=& \ds{ \frac{1}{d_{\pi} d_{\lambda}} \int_{|y|=r}  e^{-2 \langle
y, \xi \rangle} d\mu_r(y) \sum_{p,q=1}^{d_{\pi}}
\sum_{m,k=1}^{d_{\lambda}} \left| \sum_{i=1}^{d_{\pi}}
\sum_{l=1}^{d_{\lambda}} \widetilde{f^{lm}_{ip}}(\xi)
\phi_{lk}^{\lambda}(e^{-iH}) \phi_{qi}^{\pi}(e^{iH}) \right|^2 }
\end{array}\right.\eea

\vspace*{0.1in}

We have obtained an expression for one part of Lemma \ref{lem}.
Now, looking at the other part, we have
$$ f(u^{-1} \cdot x, u^{-1}k^{-1}) = \sum_{i,j=1}^{d_{\pi}}
\sum_{l,m=1}^{d_{\lambda}} f^{lm}_{ij}(x)
\phi_{lm}^{\lambda}(u^{-1}) \phi_{ji}^{\pi}(ku).$$ So, if $f$ is
holomorphic on $\mathbb{C}^n \times G,$ for $z = x+iy$ we get $$
f(e^{-iH}u^{-1} \cdot z, e^{-iH}u^{-1}k^{-1}) =
\sum_{i,j,q=1}^{d_{\pi}} \sum_{l,m=1}^{d_{\lambda}} f^{lm}_{ij}(z)
\phi_{lm}^{\lambda}(e^{-iH}u^{-1}) \phi_{jq}^{\pi}(k)
\phi_{qi}^{\pi}(ue^{iH}).$$ Again, by Schur's orthogonality
relations and similar reasoning as before, we have \beas && \int_K
\left|f(e^{-iH}u^{-1} \cdot z, e^{-iH}u^{-1}k^{-1}) \right|^2 dk \\
&=& \frac{1}{d_{\pi}} \sum_{j,q=1}^{d_{\pi}} \left|
\sum_{i=1}^{d_{\pi}} \sum_{l,m=1}^{d_{\lambda}} f^{lm}_{ij}(z)
\phi_{lm}^{\lambda}(e^{-iH}u^{-1}) \phi_{qi}^{\pi}(ue^{iH})
\right|^2 \\
&=& \frac{1}{d_{\pi}} \sum_{j,q=1}^{d_{\pi}} \left|
\sum_{i=1}^{d_{\pi}} \sum_{l,m=1}^{d_{\lambda}} f^{lm}_{ij}(z)
\phi_{lm}^{\lambda}(e^{-iH}u^{-1}) \phi_{qi}^{\pi}(e^{iH})
\right|^2. \eeas Hence, by invariance of Haar measure, we have
\beas && \int_{\mathbb{R}^n} \int_K \int_K \left|
f(e^{-iH}u^{-1} \cdot z, e^{-iH}u^{-1}k^{-1}) \right|^2 dk du dx \\
&=& \frac{1}{d_{\pi}} \sum_{j,q=1}^{d_{\pi}} \int_{\mathbb{R}^n}
\int_K \left| \sum_{i=1}^{d_{\pi}} \sum_{p,l,m=1}^{d_{\lambda}}
f^{lm}_{ij}(z) \phi_{lp}^{\lambda}(e^{-iH})
\phi_{pm}^{\lambda}(u^{-1}) \phi_{qi}^{\pi}(e^{iH}) \right|^2 du dx \\
&=& \frac{1}{d_{\pi} d_{\lambda}} \sum_{j,q=1}^{d_{\pi}}
\sum_{p,m=1}^{d_{\lambda}} \int_{\mathbb{R}^n} \left|
\sum_{i=1}^{d_{\pi}} \sum_{l=1}^{d_{\lambda}} f^{lm}_{ij}(x+iy)
\phi_{lp}^{\lambda}(e^{-iH}) \phi_{qi}^{\pi}(e^{iH}) \right|^2 dx \\
&=& \frac{1}{d_{\pi} d_{\lambda}} \sum_{j,q=1}^{d_{\pi}}
\sum_{p,m=1}^{d_{\lambda}} \int_{\mathbb{R}^n} \left|
\sum_{i=1}^{d_{\pi}} \sum_{l=1}^{d_{\lambda}}
\widetilde{f^{lm}_{ij}}(\xi) \phi_{lp}^{\lambda}(e^{-iH})
\phi_{qi}^{\pi}(e^{iH}) \right|^2 e^{-2 \langle y, \xi \rangle}
d\xi. \eeas Now by the invariance of Lebesgue measure under the
$K$-action on $\mathbb{R}^n$ we get that \beas && \int_{|y|=r}
\int_{\mathbb{R}^n} \int_K \int_K \left| f(e^{-iH}u^{-1} \cdot z,
e^{-iH}u^{-1}k^{-1}) \right|^2 dk du dx d\mu_r(y) \\
&=& \int_{|y|=r} \int_{\mathbb{R}^n} \int_K \left| f(e^{-iH} \cdot
z, e^{-iH} k) \right|^2 dk dx d\mu_r(y). \eeas Hence the lemma
follows from (\ref{eqn5}).

\end{proof}

\noindent \textit{Proof of Theorem \ref{thm}.}

To prove the theorem, it is enough to prove the orthogonality of
the components $\ds{f_{\pi}^{\lambda}(x,k) =
\sum_{i,j=1}^{d_{\pi}} \sum_{l=1}^{d_{\lambda}} f^{ll}_{ij}(x)
\overline{\phi_{ij}^{\pi}(k)}.}$ For $\pi, \lambda, \tau, \nu \in
\widehat{K},$ we have \beas && \left\langle
U_{(x+iy,ke^{iH})}^{\xi} \widehat{f_{\pi}^{\lambda}}(\xi) ,
U_{(x+iy,ke^{iH})}^{\xi}
\widehat{f_{\tau}^{\nu}}(\xi) \right\rangle_{HS} \\
&=& \sum_{\gamma \in \widehat{K}} d_{\gamma}
\sum_{p,q=1}^{d_{\gamma}} \int_K \frac{\delta_{\gamma
\pi}}{d_{\pi}} e^{i \langle x+iy, u \cdot \xi \rangle}
\sum_{i=1}^{d_{\pi}} \sum_{l,m=1}^{d_{\lambda}}
\widetilde{f^{lm}_{ip}}(\xi) \phi_{lm}^{\lambda}(e^{-iH}k^{-1}u)
\phi_{qi}^{\pi}(u^{-1}ke^{iH}) \\
&& \frac{\delta_{\gamma \tau}}{d_{\tau}} \overline{e^{i \langle
x+iy, u \cdot \xi \rangle}} \sum_{a=1}^{d_{\tau}}
\sum_{b,c=1}^{d_{\nu}} \overline{\widetilde{f^{bc}_{ap}}(\xi)
\phi_{bc}^{\nu}(e^{-iH}k^{-1}u)
\phi_{qa}^{\tau}(u^{-1}ke^{iH})} du \\
&=& 0 \txt{ if } \pi \ncong \tau. \eeas Assume $\pi \cong \tau.$
Then \beas &&  \int_{|y|=r} \left\langle U_{(x+iy,ke^{iH})}^{\xi}
\widehat{f_{\pi}^{\lambda}}(\xi) , U_{(x+iy,ke^{iH})}^{\xi}
\widehat{f_{\pi}^{\nu}}(\xi) \right\rangle_{HS} d\mu_r(y)\\
&=& \frac{1}{d_{\pi}}
\frac{J_{\frac{n}{2}-1}(2ir|\xi|)}{(2ir|\xi|)^{\frac{n}{2}-1}}
\sum_{a,i,p=1}^{d_{\pi}} \sum_{l,m=1}^{d_{\lambda}}
\sum_{b,c=1}^{d_{\nu}} \widetilde{f^{lm}_{ip}}(\xi)
\overline{\widetilde{f^{bc}_{ap}}(\xi)} \\
&& \int_K \left( \sum_{q=1}^{d_{\pi}}
\phi_{qi}^{\pi}(u^{-1}e^{iH})
\overline{\phi_{qa}^{\pi}(u^{-1}e^{iH})} \right)
\phi_{lm}^{\lambda}(e^{-iH}u) \overline{\phi_{bc}^{\nu}(e^{-iH}u)} du \\
&=& \frac{1}{d_{\pi}}
\frac{J_{\frac{n}{2}-1}(2ir|\xi|)}{(2ir|\xi|)^{\frac{n}{2}-1}}
\sum_{a,i,p,q=1}^{d_{\pi}} \sum_{l,m=1}^{d_{\lambda}}
\sum_{b,c=1}^{d_{\nu}} \widetilde{f^{lm}_{ip}}(\xi)
\overline{\widetilde{f^{bc}_{ap}}(\xi)} \phi_{qi}^{\pi}(e^{iH})
\overline{\phi_{qa}^{\pi}(e^{iH})} \\
&& \sum_{j=1}^{d_{\lambda}} \sum_{k=1}^{d_{\nu}}
\phi_{lj}^{\lambda}(e^{-iH})
\overline{\phi_{bk}^{\nu}(e^{-iH})} \int_K \phi_{jm}^{\lambda}(u) \overline{\phi_{kc}^{\nu}(u)} du \\
&=& 0 \txt{ if } \lambda \ncong \nu. \eeas On the other hand, we
have \beas && \int_K f_{\pi}^{\lambda}(e^{-iH}u^{-1} \cdot z,
e^{-iH}u^{-1}k^{-1}) \overline{f_{\tau}^{\nu}(e^{-iH}u^{-1} \cdot
z, e^{-iH}u^{-1}k^{-1})} dk \\
&=& \sum_{i,j,q=1}^{d_{\pi}} \sum_{l,m=1}^{d_{\lambda}}
\sum_{a,b,p=1}^{d_{\tau}} \sum_{s,t=1}^{d_{\nu}} f^{lm}_{ij}(z)
\overline{f^{st}_{ab}(z)} \phi_{lm}^{\lambda}(e^{-iH}u^{-1})
\overline{\phi_{st}^{\nu}(e^{-iH}u^{-1})} \\
&& \phi_{qi}^{\pi}(ue^{iH}) \overline{\phi_{pa}^{\tau}(ue^{iH})}
\int_K \phi_{jq}^{\pi}(k) \overline{\phi_{bp}^{\tau}(k)} dk\\
&=& 0 \txt{ if } \pi \ncong \tau. \eeas Assume $\pi \cong \tau.$
Then we get  \beas && \int_K \int_K
f_{\pi}^{\lambda}(e^{-iH}u^{-1} \cdot z, e^{-iH}u^{-1}k^{-1})
\overline{f_{\pi}^{\nu}(e^{-iH}u^{-1} \cdot z,
e^{-iH}u^{-1}k^{-1})} dk du \\
&=& \sum_{i,a,j=1}^{d_{\pi}} \sum_{l,m=1}^{d_{\lambda}}
\sum_{s,t=1}^{d_{\nu}} f^{lm}_{ij}(z) \overline{f^{st}_{aj}(z)}
\left( \sum_{q=1}^{d_{\pi}}  \phi_{qi}^{\pi}(e^{iH})
\overline{\phi_{pa}^{\tau}(e^{iH})} \right) \\
&& \sum_{\alpha=1}^{d_{\lambda}} \sum_{\beta=1}^{d_{\nu}}
\phi_{l\alpha}^{\lambda}(e^{-iH})
\overline{\phi_{s\beta}^{\nu}(e^{-iH})} \int_K \phi_{\alpha
m}^{\lambda}(u^{-1}) \overline{\phi_{\beta t}^{\nu}(u^{-1})} du \\
&=& 0 \txt{ if } \lambda \ncong \nu. \eeas This finishes the
proof.

\qed

It is easy to see that $$ \int_{\mathbb{R}^n} \|U_{(z,g)}^{\xi}
\widehat{f}(\xi) \|_{HS}^2 d\xi = \sum_{\sigma \in
\widehat{K_{\xi}}} d_{\sigma} \int_{\mathbb{R}^n}
\|U_{(z,g)}^{\xi, \sigma} \widehat{f}(\xi, \sigma) \|_{HS}^2
d\xi.$$ Hence we have the following corollary:

\begin{cor}
For $f \in L^2(M),$ $f$ extends holomorphically to $\mathbb{C}^n
\times G$ with $$\int_{|y| = r} \int_K \int_{\mathbb{R}^n}
|f(e^{-iH}(x+iy),e^{-iH}k)|^2 dx dk d\mu_r(y) < \infty$$ (where
$\mu_r$ is the normalized surface area measure on the sphere $\{
|y| = r \} \subseteq \mathbb{R}^n$) iff $$\ds{ \sum_{\sigma \in
\widehat{K_{\xi}}} d_{\sigma} \int_{\mathbb{R}^n} \int_{|y| = r}
\|U_{(z,g)}^{\xi, \sigma} \widehat{f}(\xi,\sigma) \|_{HS}^2
d\mu_r(y) d\xi} < \infty $$ where $z=x+iy \in \mathbb{C}^n,$ $g
\in G$ and we also have \beas && \sum_{\sigma \in
\widehat{K_{\xi}}} d_{\sigma} \int_{\mathbb{R}^n} \int_{|y| = r}
\|U_{(z,g)}^{\xi,\sigma} \widehat{f}(\xi,\sigma) \|_{HS}^2
d\mu_r(y) d\xi \\ &=& \int_{|y| = r} \int_{K} \int_{\mathbb{R}^n}
|f(e^{-iH}(x+iy),e^{-iH}k)|^2 dx dk d\mu_r(y).\eeas
\end{cor}

\textbf{Acknowledgement.} The author wishes to thank Dr. E. K.
Narayanan for his encouragement and for the many useful
discussions during the course of this work.

\end{document}